# EXPONENTIAL POLYNOMIALS AND STRATIFICATION IN THE THEORY OF ANALYTIC INEQUALITIES

BRANKO MALEŠEVIĆ[1], MILOŠ MIĆOVIĆ[1]

***Abstract.*** *This paper considers MEP - Mixed Exponential Polynomials as one class of real exponential polynomials. We introduce a method for proving the positivity of MEP inequalities over positive intervals using the Maclaurin series to approximate the exponential function precisely. Additionally, we discuss the relation between MEPs and stratified families of functions from [1] through two applications, referring to inequalities from papers [2] and [3].*

## 1. INTRODUCTION

Exponential polynomials are finite linear combinations of term $x^n e^{ax}$, where $n \in N_0 = N \cup \{0\}$ and $a \in C$ [4]. In this paper, we examine a class of real exponential polynomials defined by:

$$\sum_{k=1}^{m} \alpha_k x^{p_k} (e^{-x})^{q_k}, \tag{1}$$

where $m, k \in N$, $p_k, q_k \in N_0$ and $\alpha_k \in R$, for $x \in (a,b) \subseteq R$. These exponential polynomials we denote as MEP - Mixed Exponential Polynomials. This denomination resembles MTP - Mixed Trigonometric Polynomials from papers [5-8]. By analogy with a method for proving MTP inequalities from [5], in this paper, we introduce a method for proving MEP inequalities.

## 2. A METHOD FOR PROVING THE POSITIVITY OF MEP

In this section, we introduce the main assertions on which the method for proving the positivity of MEP over a finite interval is based. In addition, at the end of the section we discuss some properties of MEP.

Let $T_n(x)$ denotes the Maclaurin series of the exponential function $e^{-x}$:

$$T_n(x) = 1 + \sum_{k=1}^{n} \frac{(-1)^k}{k!} x^k, \tag{2}$$

where $n \in N_0$ is an order of the series and $x \in R$.

[1]University of Belgrade, School of Electrical Engineering, Department of Applied Mathematics, 11000 Belgrade, Serbia. E-mail: branko.malesevic@etf.bg.ac.rs; milos.micovic@etf.bg.ac.rs.

In the following assertions, we consider the function:

$$F_n(x) = T_n(x) - e^{-x} = 1 + \sum_{k=1}^{n} \frac{(-1)^k}{k!} x^k - e^{-x} : R \to R,$$

where $n \in N_0$ is an order of the series and $x \in R$. Then the following assertions hold:

**Lemma 1.** It holds:

$$F_n^{(k)}(x) = (-1)^k (T_{n-k}(x) - e^{-x}) \tag{3}$$

for $k \in \{0, 1, 2, \ldots, n\}$, $n \in N_0$ and $x \in R$.

**Theorem 1.** For $x > 0$ and $m \in N_0$, it holds:

$$e^{-x} < T_{2m}(x). \tag{4}$$

*Proof:* Let us consider the function $F_{2m}(x)$ for $x \in [0, \infty)$. Then $F_{2m}^{(2m)}(x) = 1 - e^{-x}$ and, since it holds:

$$F_{2m}^{(2m+1)}(x) = e^{-x} > 0,$$

we conclude that $F_{2m}^{(2m)}(x)$ is an increasing function. Hence, since $F_{2m}^{(2m)}(0) = 0$, it holds:

$$F_{2m}^{(2m)}(x) > 0$$

for $x > 0$. It follows that $F_{2m}^{(2m-1)}(x) = -1 + x + e^{-x}$ is an increasing function, and also it holds that $F_{2m}^{(2m-1)}(0) = 0$. Hence:

$$F_{2m}^{(2m-1)}(x) > 0.$$

Continuing this procedure, using Lemma 1, we can obtain that:

$$F_{2m}(x) > 0,$$

i.e.

$$e^{-x} < T_{2m}(x)$$

for $x > 0$, which proves the theorem. □

**Theorem 2.** For $x > 0$ and $m \in N$, it holds:

$$e^{-x} > T_{2m-1}(x). \tag{5}$$

*Proof:* The proof is entirely analogous to the proof of Theorem 1. □

**The description of the method.** Let us start the description of the method for proving the positivity of MEP. Let us observe MEP:

$$f(x) = \sum_{k=1}^{m} \alpha_k x^{p_k} (e^{-x})^{q_k} \tag{6}$$

on the interval $(a, b) \subseteq R^+$. First, we discuss one method for forming a polynomial function $P(x)$, so that $f(x) > P(x)$ for $x \in (a, b)$. Let us observe the addend $\alpha_k x^{p_k}(e^{-x})^{q_k}$ of the function $f(x)$. If $\alpha_k > 0$ then, for that addend, we choose a polynomial $T_{2l-1}(q_k x)$, and for such a choice it holds:

$$(e^{-x})^{q_k} > T_{2l-1}(q_k x), \tag{7}$$

where $l, q_k \in N$. If $\alpha_k < 0$ then, for that addend, we choose a polynomial $T_{2l}(q_k x)$, and for such a choice it holds:

$$(e^{-x})^{q_k} < T_{2l}(q_k x), \tag{8}$$

where $l, q_k \in N$. Thus, we determine the polynomial $P(x) = P_{l_1,\ldots,l_m}(x)$, a downward approximation of MEP, with:

$$P(x) = \sum_{k=1}^{m} \alpha_k x^{p_k} T_{\vartheta_k}(q_k x),$$

where for the degrees of the polynomial $T_{\vartheta_k}(q_k x)$, it holds:

$$\vartheta_k = \begin{cases} 2l_k - 1, & \text{for } \alpha_k > 0; \\ 2l_k, & \text{for } \alpha_k < 0. \end{cases}$$

The polynomial $P(x) = P_{l_1,\ldots,l_m}(x)$, where $l_1, \ldots, l_m \in N$ and $m \in N$, by the method of its formation, satisfies:

$$f(x) = \sum_{k=1}^{m} \alpha_k x^{p_k}(e^{-x})^{q_k} > P(x) = \sum_{k=1}^{m} \alpha_k x^{p_k} T_{\vartheta_k}(q_k x).$$

By MEP inequality, we mean an inequality of the form:

$$f(x) > 0$$

for $x \in (a, b)$. The method of proof of MEP inequality is based on the following assertion:

**Theorem 3.** If there exist indices $l_1, \ldots, l_m$, so that for the polynomial $P(x) = P_{l_1,\ldots,l_m}(x)$, it holds:

$$P(x) > 0$$

for $x \in (a, b)$, then:

$$f(x) > 0$$

for $x \in (a, b)$.

**Corollary 1.** According to the previous assertion, the proof of MEP inequality $f(x) > 0$ reduces to the proof of a polynomial inequality $P(x) > 0$, if such an inequality is possible for some choice of indices $l_1, \ldots, l_m$, for a polynomial with rational coefficients. In that case, the proof of an inequality P(x) > 0 is algorithmically decidable using the Sturm's theorem, assuming that $(a, b)$ is an interval with rational ends [9] (see Theorem 4.2 b).

**Remark 4.** The function $f(x) = \sum_{k=1}^{m} \alpha_k x^{p_k} (e^{-x})^{q_k}$ can be considered in the case when $q_k < 0$. These exponential polynomials are not MEPs, but inequalities $f(x) > 0$, which we denote as exponential polynomial inequalities, can be reduced to MEP inequalities by multiplying the inequality $f(x) > 0$ by $e^{q_k x}$ (multiple times) whenever $q_k < 0$. With this procedure, every exponential polynomial inequality reduces to MEP inequality.

**Remark 5.** Let us note that the presented method can be improved in some cases. If, in the notation of MEP, there are multiple addends with the same exponent $(e^{-x})^{q_k}$ that form a polynomial, and if such a polynomial is positive or negative over the interval $(a, b)$, it is not necessary to perform approximations of type (7) and (8) for each addend. In that case, it is sufficient to consider the product of such a polynomial and $(e^{-x})^{q_k}$ function in order to determine a downward/upward polynomial approximation.

At the end of this section, we present the assertion about the order of the Maclaurin approximations of the function $e^{-x}$, which can be applied in proving MEP inequalities.

**Theorem 6.** Let $T_i(x)$ be the sequence of Maclaurin polynomials of the function $e^{-x}$ ($i \in N_0$) which we consider for $x \in [0, \infty)$. It holds:

**(i)** For each Maclaurin polynomial $T_{2m-1}(x)$ of odd degree, there exists exactly one root $c_{2m-1} \in (0, \infty)$ ($m \in N$). Moreover, the following order: $c_1 = 1 < c_3 < \cdots$ is true.

**(ii)** Over the interval $(0, 1)$, the following inequalities hold:

$$T_1(x) < T_3(x) < \cdots < T_{2m-1}(x) < e^{-x} < T_{2m}(x) < \cdots < T_2(x) < T_0(x).$$

**(iii)** The following equality holds:

$$T_n(x) = T_{n+2}(x) \quad \text{iff} \quad x = 0 \ \lor \ x = n + 2,$$

as well as:

$$T_{2m-1}(x) < T_{2m+1}(x), for \ x \in (0, d) \quad \text{and} \quad T_{2m-1}(x) > T_{2m+1}(x), for \ x \in (d, +\infty),$$

where $d = 2m + 1$, and:

$$T_{2m-2}(x) > T_{2m}(x), for \ x \in (0, d) \quad \text{and} \quad T_{2m-2}(x) < T_{2m}(x), for \ x \in (d, +\infty),$$

where $d = 2m$.

*Proof:* **(i)** Let us first examine the monotonicity of the function $T_{2m+1}(x)$. The first derivative of the function $T_{2m+1}(x)$ is:

$$T'_{2m+1}(x) = -1 + x - \frac{1}{2!}x^2 + \cdots - \frac{1}{(2m)!}x^{2m} = -T_{2m}(x).$$

Based on Theorem 1, it follows that: $T_{2m}(x) > e^{-x} > 0$. Hence, we conclude that $T'_{2m+1}(x) < 0$, and therefore $T_{2m+1}(x)$ is a decreasing function.

For the polynomial $T_{2m+1}(x) = -\frac{1}{(2m+1)!}x^{2m+1} + \cdots + \frac{1}{2!}x^2 - x + 1$, it holds:

$$T_{2m+1}(0) = 1 \quad \text{and} \quad \lim_{x \to \infty} T_{2m+1}(x) = -\infty.$$

Based on the continuity of the polynomial function $T_{2m+1}(x)$, we conclude that there exists exactly one root $c_{2m+1} \in (0, +\infty)$. It is evident that $c_1 = 1$. It holds:

$$T_{2m-1}(x) < T_{2m+1}(x) \Leftrightarrow x \in (0, 2m+1).$$

Moreover,

$$T_{2m+1}(2m+1) = 1 - (2m+1) + \frac{(2m+1)^2}{2!} - \frac{(2m+1)^3}{3!} + \cdots + \frac{(2m+1)^{2m}}{(2m)!} - \frac{(2m+1)^{2m+1}}{(2m+1)!} < 0$$

since every two by two consecutive addends are non-positive. From this, we can conclude that:

$$c_{2m-1} < c_{2m+1}$$

for each $m \in N$. Thus $c_1 = 1 < c_3 < \cdots$.

**(ii)** Based on the equivalences: $T_{2m-1}(x) < T_{2m+1}(x) \Leftrightarrow x \in (0, 2m+1)$ and $T_{2m-2}(x) > T_{2m}(x) \Leftrightarrow x \in (0, 2m)$ (which are evident), and Theorems 1 and 2, the assertion follows.

**(iii)** The stated equality and inequalities are evident. □

## 3. STRATIFICATION AND MEP INEQUALITIES

In this section, we first summarize the basic facts about stratified families of functions according to the paper [1]. The family of functions $\varphi_p(x)$, where $x \in (a,b) \subseteq R^+$, $p \in R^+$, is increasingly stratified with respect to parameter $p$ if:

$$(\forall p_1, p_2 \in R^+)\ p_1 < p_2 \Leftrightarrow \varphi_{p_1}(x) < \varphi_{p_2}(x)$$

holds for each $x \in (a,b)$, and conversely, decreasingly stratified with respect to parameter p if:

$$(\forall p_1, p_2 \in R^+)\ p_1 < p_2 \Leftrightarrow \varphi_{p_1}(x) > \varphi_{p_2}(x)$$

holds for each $x \in (a,b)$. Note that in the paper [1], it is considered, for stratified families of functions $\varphi_p(x)$, when it is possible to determine the unique value $p \in R^+$, for which the infimum of the error $d(p) = sup_{x \in (a,b)}|\varphi_p(x)|$ is attained. For such a value $p = p_0$, the function $\varphi_p(x)$ is called *the minimax approximant* on $(a,b)$ [1].

Furthermore, we will present, through two applications, the application of MEP and stratification on examples of inequalities from papers [2] and [3].

*3.1. APPLICATION 1*

This subsection provides the application of the previously presented method for proving the positivity of MEP on the example of improving the results from the paper Y.-D. Wu and Z.-H. Zhang [2].

We quote the result of that paper.

**Theorem 7.** Let $0 < x < 1$. Then the following inequality holds:

$$\frac{1}{x^2} - \frac{1}{12} < \frac{e^{-x}}{(1-e^{-x})^2} \qquad (9)$$

and the constant $1/12$ is the best possible.

Naturally, regarding the previous inequality, we take into consideration the following family of functions:

$$\varphi_p(x) = \frac{1}{x^2} - \frac{e^{-x}}{(1-e^{-x})^2} - p \qquad (10)$$

for $x \in (0, 1)$ and $p > 0$. Let us note that exist:

$$\varphi_p(0+) = \frac{1}{12} - p \in R \quad \text{and} \quad \varphi_p(1) = 1 - \frac{e^{-1}}{(1-e^{-1})^2} - p \in R.$$

The following assertions hold, which provide improvements to Theorem 7.

**Lemma 2.** The family of functions:

$$\varphi_p(x) = \frac{1}{x^2} - \frac{e^{-x}}{(1-e^{-x})^2} - p \qquad (for\ x \in (0,1))$$

is decreasingly stratified with respect to parameter $p \in R^+$.

*Proof:* $\frac{\partial \varphi_p(x)}{\partial p} = -1 < 0$. □

**Statement 1.** Let

$$A = \frac{e^2 - 3e + 1}{e^2 - 2e + 1} = 0.079326\ldots \quad and \quad B = \frac{1}{12} = 0.083333\ldots.$$

(*i*) If $p \in (0, A]$, then:

$$x \in (0,1) \Rightarrow \frac{e^{-x}}{(1-e^{-x})^2} < \frac{1}{x^2} - A \leq \frac{1}{x^2} - p.$$

(*ii*) If $p \in (A, B)$, then the equality:

$$\varphi_p(x) = \frac{1}{x^2} - \frac{e^{-x}}{(1-e^{-x})^2} - p = 0$$

has the unique solution $x_0^{(p)}$ and it holds:

$$x \in \left(0, x_0^{(p)}\right) \Rightarrow \frac{e^{-x}}{(1-e^{-x})^2} < \frac{1}{x^2} - p$$

and:

$$x \in \left(x_0^{(p)}, 1\right) \Rightarrow \frac{e^{-x}}{(1-e^{-x})^2} > \frac{1}{x^2} - p.$$

(*iii*) If $p \in [B, \infty)$, then:

$$x \in (0, 1) \Rightarrow \frac{e^{-x}}{(1-e^{-x})^2} > \frac{1}{x^2} - B \geq \frac{1}{x^2} - p.$$

(*iv*) There is exactly one solution to the equation:

$$\varphi_p(0+) = \varphi_p(1-)$$

with respect to parameter $p \in (A, B)$, determined with:

$$p_0 = \frac{13e^2 - 38e + 13}{24e^2 - 48e + 24} = 0.081329 \ldots.$$

For value:

$$d_0 = \varphi_p(0+) = |\varphi_{p_0}(1)| = \frac{-11e^2 + 34e - 11}{24e^2 - 48e + 24} = 0.0020034 \ldots$$

it holds:

$$d_0 = \inf_{p \in R^+} \sup_{x \in (0,1)} |\varphi_p(x)|.$$

(*v*) For value $p_0 = \frac{13e^2 - 38e + 13}{24e^2 - 48e + 24} = 0.081329 \ldots$ the minimax approximant of the family of functions is:

$$\varphi_{p_0}(x) = \frac{1}{x^2} - \frac{e^{-x}}{(1-e^{-x})^2} - \frac{13e^2 - 38e + 13}{24e^2 - 48e + 24},$$

which determine the (minimax) approximation:

$$\frac{e^{-x}}{(1-e^{-x})^2} \approx \frac{1}{x^2} - 0.081329 \,.$$

*Proof:* Let $p > 0$. Let us examine the monotonicity of the function:

$$\varphi_p(x) = \frac{1}{x^2} - \frac{e^{-x}}{(1-e^{-x})^2} - p$$

on the interval $(0, 1)$. The first derivative of the function $\varphi_p(x)$ is:

$$\frac{\partial \varphi_p(x)}{\partial x} = \frac{-2}{x^3} + \frac{e^{-x}}{(1-e^{-x})^2} + \frac{2(e^{-x})^2}{(1-e^{-x})^3}$$

$$= \frac{-e^{-2x}x^3 - e^{-x}x^3 - 2e^{-3x} + 6e^{-2x} - 6e^{-x} + 2}{x^3(-1+e^{-x})^3}.$$

Let us observe $g(x) = -e^{-2x}x^3 - e^{-x}x^3 - 2e^{-3x} + 6e^{-2x} - 6e^{-x} + 2$ as one MEP. Based on this MEP, we form two simpler MEPs:

$$g_1(x) = 6e^{-2x} + 2$$

and

$$g_2(x) = e^{-2x}x^3 + e^{-x}x^3 + 2e^{-3x} + 6e^{-x}.$$

Then:

$$g(x) = g_1(x) - g_2(x),$$

where the polynomial $g_1(x)$ consists of those addends of MEP $g(x)$ for which $\alpha_k > 0$, while the polynomial $g_2(x)$ consists of those addends of $g(x)$ for which $\alpha_k < 0$.

Let us approximate MEPs $g_1(x)$ and $g_2(x)$ by Maclaurin polynomials. If we approximate $e^{-2x}$ by the Maclaurin polynomial of degree 9, then $g_1(x) = 6e^{-2x} + 2$ has the downward polynomial approximation:

$$G_1(x) = -\frac{8}{945}x^9 + \frac{4}{105}x^8 - \frac{16}{105}x^7 + \frac{8}{15}x^6 - \frac{8}{5}x^5 + 4x^4 - 8x^3 + 12x^2 - 12x + 8$$

for which $g_1(x) > G_1(x)$ for $x \in (0, 1)$ based on Theorem 2. Similarly, if we approximate $e^{-x}$, $e^{-2x}$ and $e^{-3x}$ by the Maclaurin polynomial of degree 12, then $g_2(x) = e^{-2x}x^3 + e^{-x}x^3 + 2e^{-3x} + 6e^{-x}$ has the upward polynomial approximation:

$$G_2(x) = \frac{4097}{479001600}x^{15} - \frac{683}{13305600}x^{14} + \frac{41}{145152}x^{13} + \frac{287}{356400}x^{12}$$

$$- \frac{3329}{1330560}x^{11} + \frac{1051}{151200}x^{10} - \frac{551}{30240}x^9 + \frac{17}{336}x^8 - \frac{9}{56}x^7$$

$$+ \frac{8}{15}x^6 - \frac{8}{5}x^5 + 4x^4 - 8x^3 + 12x^2 - 12x + 8$$

for which $g_2(x) < G_2(x)$ for $x \in (0, 1)$ based on Theorem 1.

Let:

$$G(x) = G_1(x) - G_2(x)$$

$$= -\frac{4097}{479001600}x^{15} + \frac{683}{13305600}x^{14} - \frac{41}{145152}x^{13} - \frac{287}{356400}x^{12}$$

$$+ \frac{3329}{1330560}x^{11} - \frac{1051}{151200}x^{10} + \frac{59}{6048}x^9 - \frac{1}{80}x^8 + \frac{1}{120}x^7$$

It holds that $g(x) > G(x)$ for $x \in (0,1)$, since $g_1(x) > G_1(x)$ and $g_2(x) < G_2(x)$ on the observed interval $(0,1)$. Based on the Sturm's theorem, it can be concluded that the polynomial function $G(x)$ does not have any root on the interval $(0,1)$. It holds: $G(1/2) = \frac{463573639}{15695924428800} = 0.000029534\ldots > 0$. Based on that, it holds: $G(x) > 0$ for $x \in (0,1)$, and therefore, according to Theorem 3, for MEP $g(x)$ it is proved:

$$g(x) > 0$$

on the interval $(0,1)$.

Since $g(x) > 0$, we conclude that:

$$\frac{\partial \varphi_p(x)}{\partial x} < 0,$$

and therefore $\varphi_p(x)$ is a decreasing function on the interval $(0,1)$ for each value of the parameter $p$.

Since $\varphi_p(x)$ is a decreasing function on the interval $(0,1)$ and since $\varphi_A(1-) = 0$ and $\varphi_B(0+) = 0$, the assertions $(i)$, $(ii)$ and $(iii)$ are proved.
Due to the fact that $\varphi_p(x)$, for each value of the parameter $p$, is a decreasing function on the interval $(0,1)$ and that the family of functions $\varphi_p(x)$ is stratified, the minimax approximant of the family of functions is obtained for value $p = p_0$ for which it holds:

$$\varphi_{p_0}(0+) = -\varphi_{p_0}(1-).$$

Solving this equation using the Computer Algebra System *Maple* obtains the stated value:

$$p_0 = \frac{13e^2 - 38e + 13}{24e^2 - 48e + 24} = 0.081329\ldots$$

for which it holds:

$$d_0 = inf_{p \in R^+} sup_{x \in (0,1)} |\varphi_p(x)|.$$

This finishes the proof. □

**Remark 8.** MEPs $g_1(x)$ and $g_2(x)$ could have been approximated, in the proof, by Maclaurin polynomials of another degree. The best choice, which provides the minimal degree of the polynomial $G(x)$, would be to approximate the addend $6e^{-2x}$ of the polynomial $g_1(x)$ by the Maclaurin polynomial of degree 11, while, regarding the polynomial $g_2(x)$, it would be the best to approximate the addends $e^{-2x}x^3$, $e^{-x}x^3$, $2e^{-3x}$, $6e^{-x}$ by the Maclaurin polynomials of degrees 8, 6, 12, and 8, respectively. With this particular selection of Maclaurin approximations, the degree of the polynomial $G(x)$ would be 12.

**Corollary 2.** Let $A = \frac{e^2 - 3e + 1}{e^2 - 2e + 1} = 0.079326 \ldots$ and $B = \frac{1}{12} = 0.083333 \ldots$. Then for each $0 < x < 1$, it holds:

$$\frac{1}{x^2} - p_2 \leq \frac{1}{x^2} - B < \frac{e^{-x}}{(1 - e^{-x})^2} < \frac{1}{x^2} - A \leq \frac{1}{x^2} - p_1, \qquad (11)$$

where $p_1 \leq A$ and $p_2 \geq B$.

**Corollary 3.** For $x \in (0, 1)$, the following approximation holds:

$$\frac{e^{-x}}{(1 - e^{-x})^2} \approx \frac{1}{x^2} - M, \qquad (12)$$

where:

$$M = \frac{13e^2 - 38e + 13}{24e^2 - 48e + 24} = 0.081329 \ldots.$$

The previous approximation (12), obtained based on the minimax approximant $\varphi_{p_0}(x)$, is considered to be the minimax approximation.

### 3.2. APPLICATION 2

In this subsection, we present a new and simple proof of the inequality, which reduces to the MEP inequality, from the paper C. Chesneau, Y. J. Bagul and R. M. Dhaigude [3], based on the properties of stratified families of functions. Additionally, an extension of the claim is obtained.

If we denote function *sign* in a standard manner: $sing(a) = -1$ if $a < 0$, $sing(a) = 0$ if $a = 0$ and $sing(a) = 1$ if $a > 0$, then the following statement is true [3]:

**Theorem 9.** Let $x \in [0, 1]$ and $a \in R$. Then the following inequality holds:

$$sing(a)e^{ax} \leq sing(a)(ax(1 - x) + x^2(e^a - 1) + 1). \qquad (13)$$

*Proof:* We need to consider the following three cases: $a < 0$, $a = 0$ and $a > 0$.

**(1)** The case $a < 0$ reduces to proving the inequality:

$$-e^{ax} \leq -ax(1 - x) - x^2(e^a - 1) - 1$$

on the interval $[0, 1]$. Let us introduce the substitution $\alpha = -a$. The previous inequality reduces to the inequality:

$$-e^{-\alpha x} \leq \alpha x(1 - x) - x^2(e^{-\alpha} - 1) - 1$$

for $\alpha > 0$ on the interval $[0, 1]$. Let us consider the family of functions:

$$\varphi_\alpha(x) = e^{-\alpha x} + \alpha x(1 - x) - x^2(e^{-\alpha} - 1) - 1$$

on the interval $(0, 1)$ for $\alpha \geq 0$. Each of these functions $\varphi_\alpha(x)$ is MEP. The first derivative of the function $\varphi_\alpha(x)$ with respect to $\alpha$ is MEP:

$$\frac{\partial \varphi_\alpha(x)}{\partial \alpha} = -xe^{-\alpha x} + x(1 - x) + x^2 e^{-\alpha}.$$

The second derivative of the function $\varphi_\alpha(x)$ with respect to $\alpha$ is also MEP:

$$\frac{\partial^2 \varphi_\alpha(x)}{\partial \alpha^2} = x^2(e^{-\alpha x} - e^{-\alpha}).$$

It holds: $\frac{\partial^2 \varphi_\alpha(x)}{\partial \alpha^2} > 0$ for $x \in (0, 1)$. Therefore, the family of functions $\frac{\partial \varphi_\alpha(x)}{\partial \alpha}$ is increasingly stratified on the interval $(0, 1)$. Since $\frac{\partial \varphi_\alpha(x)}{\partial \alpha} = -xe^{-\alpha x} + x(1 - x) + x^2 e^{-\alpha} = 0$ for $\alpha = 0$, based on the definition of the increasingly stratified families of functions, for $\alpha > 0$, it holds:

$$\frac{\partial \varphi_\alpha(x)}{\partial \alpha} > 0$$

on the interval $(0, 1)$. Hence, the family of functions $\varphi_\alpha(x)$ is increasingly stratified on the interval $(0, 1)$. Again, considering that $\varphi_\alpha(x) = e^{-\alpha x} + \alpha x(1 - x) - x^2(e^{-\alpha} - 1) - 1 = 0$ for $\alpha = 0$ and based on the definition of stratification, it holds:

$$\varphi_\alpha(x) > 0$$

on the interval $(0, 1)$ for $\alpha > 0$. The equality is evident for $x = 0$ and $x = 1$, which completes the proof for the case $a < 0$.

**(2)** For $a = 0$ the inequality is evident.

**(3)** The case $a > 0$ reduces to proving the inequality:

$$e^{ax} \leq ax(1 - x) + x^2(e^a - 1) + 1$$

on the interval $[0, 1]$. This inequality can be reduced to MEP inequality by multiplying the inequality by $e^{-ax}$, but it can be proved analogously to the proof of case (1) even without reducing it to the MEP inequality. □

In the following statement, we present an extension of Theorem 9.

**Statement 2.** Let $a \in R$. Then, for $x \in (-\infty, 1]$ it holds:

$$sing(a)e^{ax} \leq sing(a)(ax(1 - x) + x^2(e^a - 1) + 1), \qquad (14)$$

and for $x \in (1, +\infty)$ it holds:

$$sing(a)e^{ax} \geq sing(a)(ax(1 - x) + x^2(e^a - 1) + 1). \qquad (15)$$

*Proof:* The proof is entirely analogous to the proof of Theorem 9. □

## 4. CONCLUSION

In the Theory of Analytic Inequalities, there are numerous examples of inequalities with an exponential function [10], see also [11-14]. In recently published papers, authors researched exponential polynomials and related inequalities [3, 15-20]. This paper describes one class of exponential polynomials - MEPs, and stratification in the theory of analytic inequalities. MEP inequalities are significant since any exponential polynomial inequality can be reduced to MEP inequality. This paper introduces one method for proving the positivity of MEP over a positive interval. Furthermore, the application of this method is shown through the proof of the inequality from the paper [2]. A simple proof of the inequality from the paper [3], which reduces to the MEP inequality, is given using the properties of stratified families of functions. New results were obtained regarding both inequalities. Since the described method for proving the positivity of MEP is widely universal, it can be also used for proving other MEPs. Proving inequalities based on the properties of stratified families of functions can be applied to some MEP inequalities, as well as to some inequalities that are not MEP.

***Acknowledgments***: *The authors are supported in part by the Serbian Ministry of Education, Science and Technological Development, under project 451-03-47/2023-01/200103.*

## REFERENCES

[1] Malešević, B., Mihailović, B., *Applicable Analysis and Discrete Mathematics*, **15**(2), 486, 2021.
[2] Wu, Y. D., Zhang, Z. H., *The Best Constant for an Inequality*, Victoria University Research Repository, Collection **7**(1), 2004. https://vuir.vu.edu.au/id/eprint/18041
[3] Chesneau, C., Bagul, Y. J., Dhaigude, R. M., *Asia Pacific Journal of Mathematics*, **9**(6), 1, 2022.
[4] Laird, P. G., *Journal of the Australian Mathematical Society*, **17**(3), 257, 1972.
[5] Malešević, B., Makragić, M., *Journal of Mathematical Inequalities*, **10**(3), 849, 2016.
[6] Malešević, B., Banjac, B., Šešum-Čavić, V., Korolija, N., *One algorithm for testing annulling of mixed trigonometric polynomial functions on boundary points*, Proceedings of 30th TELFOR Conference, 1, 2022.
[7] Malešević, B., Banjac, B., *Automated Proving Mixed Trigonometric Polynomial Inequalities*, Proceedings of 27th TELFOR conference, 1, 2019.
[8] Lutovac, T., Malešević, B., Mortici, C., *Journal of Inequalities and Applications*, **2017**(116), 1, 2017.
[9] Cutland, N., *Computalibity: An introduction to recursive funtion theory*, Cambridge University Press, Bath, 109, 1980.
[10] Mitrinović, D. S., *Analytic Inequalities*, Springler-Verlag, Berlin, 1970.
[11] Bagul, Y. J., *Journal of Mathematical Inequalities*, **11**(3), 695, 2017.
[12] Bagul, Y. J., Chesneau, C., *CUBO A Mathematical Journal*, **21**(1), 21, 2019.
[13] Malešević, B., Lutovac, T., Banjac, B., *Filomat*, **32**(20), 6921, 2018.
[14] Milovanović, G. V., Rassias, M. (Eds.), *Analytic Number Theory, Approximation Theory and Special Functions*, Springer, New York, 2014.
[15] Chesneau, C., *Jordan Journal of Mathematics and Statistics*, **11**(3), 273, 2018.
[16] Chesneau, C., Navarro, F., *Journal of Mathematical Modeling*, **7**(2), 221, 2019.
[17] Bagul, Y. J., Dhaigude, R. M., Kostić, M., Chesneau, C., *Axioms*, **10**(4), 308, 2021.

[18] Bagul, Y. J., Chesneau, C., Dhaigude, R. M., *Mathematical Analysis and its Contemporary Applications*, **5**(1), 85, 2023.
[19] Qi, F., *Mathematical Inequalities & Application*, **23**(1), 123, 2020.
[20] Bae, J., *Mathematical Inequalities & Applications*, **16**(3), 763, 2013.